\documentclass[11pt]{amsproc}%
\usepackage{amssymb}
\usepackage{amsfonts}
\usepackage{amsmath}
\usepackage{graphicx}%
\providecommand{\U}[1]{\protect\rule{.1in}{.1in}}
\theoremstyle{plain}

\newtheorem{lemma}{Lemma}

\newtheorem{theorem}{Theorem}
\numberwithin{equation}{section}
\begin{document}
\title[Approximation by analytic maps]{Approximation of functions and their derivatives by analytic maps on certain
Banach spaces}
\author{D. Azagra}
\address{ICMAT (CSIC-UAM-UC3-UCM), Departamento de An{{\'a}}lisis Matem{{\'a}}tico,
Facultad Ciencias Matem{{\'a}}ticas, Universidad Complutense, 28040, Madrid, Spain}
\email{daniel\_azagra@mat.ucm.es}
\urladdr{}
\author{R. Fry}
\curraddr{Department of Mathematics and Statistics, Thompson Rivers University,
Kamloops, B.C., Canada}
\email{rfry@tru.ca}
\urladdr{}
\author{L. Keener}
\address{Department of Mathematics and Statistics, University of Northern British
Columbia, Prince George, B.C., Canada}
\email{keener@unbc.ca}
\urladdr{}
\thanks{}
\thanks{The second named author partly supported by NSERC (Canada).}
\thanks{}
\date{}
\dedicatory{ }
\begin{abstract}
Let $X$ be a separable Banach space which admits a separating polynomial; in
particular $X$ a Hilbert space. Let $f:X\rightarrow\mathbb{R}$ be bounded,
Lipschitz, and $C^{1}$ with uniformly continuous derivative. Then for each
$\varepsilon>0,$ there exists an analytic function $g:X\rightarrow\mathbb{R}$
with $\left\vert g-f\right\vert <\varepsilon$ and $\left\Vert g^{\prime
}-f^{\prime}\right\Vert <\varepsilon.$

\end{abstract}
\maketitle

\section{Introduction}

\medskip

The problem of approximating a smooth function and its derivatives by a
function of higher order smoothness on a Banach space $X$ has been
investigated by several authors, although the number of such results is
limited. When $X$ is finite dimensional excellent results are known, and in
fact Whitney in his classical paper \cite{W} provides essentially a complete
answer by showing: for every $C^{k}$ function $f:\mathbb{R}^{n}\rightarrow
\mathbb{R}^{m}$ and every continuous $\varepsilon:\mathbb{R}^{n}%
\rightarrow(0,+\infty)$ there exists a real analytic function $g$ such that
$\Vert D^{j}g(x)-D^{j}f(x)\Vert\leq\varepsilon(x)$ for all $x\in\mathbb{R}%
^{n}$ and $j=1,...,k$. This is the so-called $C^{k}$ fine approximation of $f$.

\medskip

The first results for $X$ infinite dimensional concern the smooth,
non-analytic case, and are due to Moulis \cite{M}. She proves, in particular,
a $C^{1}$ fine approximation theorem; namely, that for $X=c_{0}\ $or $l_{p}$
with $1<p<\infty,$ and $Y$ an arbitrary Banach space, given a $C^{1}$ map
$f:X\rightarrow Y,$ and a continuous function $\varepsilon:X\rightarrow\left(
0,\infty\right)  ,$ there exists a $C^{k}$ smooth map $g:X\rightarrow Y$
(where the optimal value of $k\geq1$ depends on the choice of $X$) such that
$\left\vert g\left(  x\right)  -f\left(  x\right)  \right\vert <\varepsilon
\left(  x\right)  $ and $\left\Vert g^{\prime}\left(  x\right)  -f^{\prime
}\left(  x\right)  \right\Vert <\varepsilon\left(  x\right)  .$ This result
was later extended in \cite{AFGJL} to the case where $X$ has an unconditional
basis and admits a Lipschitz, $C^{k}$ smooth bump function. Further work along
this line can be found in \cite{HJ}, where using ideas from \cite{F1} it is
shown that for certain range spaces $Y,$ one can relax the conditions on $X$
in \cite{AFGJL} and, for example, take $X$ to be merely separable.

\medskip

It is important to note that all the results mentioned above require, in a
very essential way, a theorem concerning the approximation of Lipschitz
functions $f$ by more regular, Lipschitz functions $g$, where the Lipschitz
constant of $g$ is fixedly proportional to the Lipschitz constant of $f,$
regardless of the precision in the approximation. In \cite{M} and \cite{AFGJL}
this is achieved by reducing the problem to the finite dimensional case using
the unconditional basis, but otherwise without this reduction traditional
methods of smooth approximation, such as smooth partitions of unity, do not
work in addressing this problem. A new approach was found in \cite{F1}, and
further developed in \cite{AFM}, \cite{F2}, \cite{AFK2}, and \cite{HJ}. This
technique has been called the method of \textit{sup-partitions of unity} in
\cite{HJ}. It seems that $C^{k}$ fine approximation must rely on such results.

\medskip

Concerning $C^{k}$ fine approximation by analytic functions for $X$ infinite
dimensional, nothing is known. In view of the remarks in the preceding
paragraph, it would appear that first one needs the ability to approximate
Lipschitz functions by Lipschitz, analytic functions with good control over
the Lipschitz constant. That is, one requires a kind of analytic sup-partition
of unity. Only very recently has this been possible with the work of
\cite{AFK1}, where it is proven that if $X$ is separable and admits a
separating polynomial, then for every Lipschitz function $f:X\rightarrow
\mathbb{R}$ and $\varepsilon>0$ there exists a Lipschitz, analytic function
$g:X\rightarrow\mathbb{R}$ with $\left\vert f-g\right\vert <\varepsilon$ and
Lip$\left(  g\right)  \leq C$Lip$\left(  f\right)  ,$ where the constant $C>1$
depends only on $X$ (for a precursor to this work see \cite{FK}). Using this,
we are able in this note to give the first results on the $C^{1}$ fine
analytic approximation problem in infinite dimensions. We remark that this
work is new even for $X$ a separable Hilbert space. We establish,

\begin{theorem}
\label{Main Theorem}Let $X$ be a separable Banach space which admits a
separating polynomial. Let $f:X\rightarrow\mathbb{R}$ be bounded and
Lipschitz, with uniformly continuous derivative, and $\varepsilon>0.$ Then
there exists an analytic function $g:X\rightarrow\mathbb{R}$ such that
$\left\vert f-g\right\vert <\varepsilon$ and $\left\Vert f^{\prime}-g^{\prime
}\right\Vert <\varepsilon.$
\end{theorem}

\medskip

\noindent Our notation is standard, with $X$ denoting a Banach space, and an
open ball with centre $x$ and radius $r$ denoted $B_{r}(x).$ If $\left\{
f_{j}\right\}  _{j}$ is a sequence of Lipschitz functions defined on $X,$ then
we will at times say this family is \textit{equi-Lipschitz} if there is a
common Lipschitz constant for all $j.$ A \textit{homogeneous polynomial of
degree }$n$ is a map,\textit{\ }$P:X\rightarrow\mathbb{R},$ of the form
$P\left(  x\right)  =A\left(  x,x,...,x\right)  ,$ where $A:X^{n}%
\rightarrow\mathbb{R}$ is $n-$multilinear and continuous. For $n=0$ we take
$P$ to be constant. A \textit{polynomial of degree }$n$ is a sum $\sum
_{i=0}^{n}P_{i}\left(  x\right)  ,$ where $i\geq1$ the $P_{i}$ are
$i$-homogeneous polynomials.

\medskip

\noindent Let $X$ be a Banach space, and $G\subset X$ an open subset. A
function $f:G\rightarrow\mathbb{R}$ is called \textit{analytic }if for every
$x\in G,$ there are a neighbourhood $N_{x},$ and homogeneous polynomials
$P_{n}^{x}:X\rightarrow\mathbb{R}$ of degree $n$, such that
\[
f\left(  x+h\right)  =\sum_{n\geq0}P_{n}^{x}\left(  h\right)
\;\text{provided\ }x+h\in N_{x}.
\]
\noindent Further information on polynomials may be found, for example, in
\cite{SS}$.$

For a Banach space $X,$ we define its (Taylor) complexification $\widetilde
{X}=X\bigoplus iX$ with norm
\[
\left\Vert x+iy\right\Vert _{\widetilde{X}}=\sup_{0\leq\theta\leq2\pi
}\left\Vert \cos\theta\ x-\sin\theta\ y\right\Vert _{X}=\sup_{T\in X^{\ast
},\Vert T\Vert\leq1}\sqrt{T(x)^{2}+T(y)^{2}}.
\]
If $L:E\rightarrow F$ is a continuous linear mapping between two real Banach
spaces then there is a unique continuous linear extension $\widetilde
{L}:\widetilde{E}\rightarrow\widetilde{F}$ of $L$ (defined by $\widetilde
{L}(x+iy)=L(x)+iL(y)$) such that $\Vert\widetilde{L}\Vert=\Vert L\Vert$. For a
continuous $k$-homogeneous polynomial $P:E\rightarrow\mathbb{R}$ there is also
a unique continuous $k$-homogeneous polynomial $\widetilde{P}:\widetilde
{E}\rightarrow\mathbb{C}$ such that $\widetilde{P}=P$ on $E\subset
\widetilde{E}$ (but the norm of $P$ is not generally preserved: one has that
$\Vert\widetilde{P}\Vert\leq2^{k-1}\Vert P\Vert$). It follows that if
$q\left(  x\right)  $ is a continuous polynomial on $X,$ there is a unique
continuous polynomial $\widetilde{q}\left(  z\right)  =\widetilde{q}\left(
x+iy\right)  $ on $\widetilde{X}$ where for $y=0$ we have $\widetilde{q}=q.$
For more information on complexifications (and polynomials) we recommend
\cite{MST}. In the sequel, all extensions of functions from $X$ to
$\widetilde{X},$ as well as subsets of $\widetilde{X},$ will be embellished
with a tilde.

\section{Main Results}

\subsection{The functions $\varphi_{n}$}

To prove Theorem 1, we start with a lemma which is a variation of \cite[Lemma
3]{AFK1}, where here we have made three changes: added part $\left(
4^{\prime}\right)  ;$ included constants $M_{n}$ for the estimate in $\left(
5\right)  ;$ and relaxed the condition that $r\geq1$ to $r>0.$ To obtain
$\left(  4^{\prime}\right)  ,$ we replace the function $b_{n}$ in the proof of
\cite[Lemma 3]{AFK1} with a $C^{1}$ version; the change in $\left(  5\right)
$ is easily handled; and requiring merely $r>0$ means that certain constants
will depend on $r,$ but as we shall apply the lemma with $r$ fixed throughout,
this causes no problem.

\medskip

\noindent First we need some definitions and notation. If $X$ posseses an
$n^{\text{th}}$ order separating polynomial, then it admits a $2n$-homogeneous
polynomial $q$ such that
\begin{equation}
\left\Vert x\right\Vert ^{2n}\leq q\left(  x\right)  \leq A\left\Vert
x\right\Vert ^{2n},
\end{equation}
for some $A>1$ (see e.g., \cite{AFK1}). In \cite[Lemma 2]{AFK1} it is proved
that the function $Q\left(  x\right)  =\left(  q\left(  x\right)  +1\right)
^{1/2n}-1$ satisfies:

\begin{enumerate}
\item $Q$ is (real) analytic on $X$,

\item $Q$ is Lipschitz on $X$, where we can take Lip$\left(  Q\right)  >1,$

\item $\inf\left\{  Q\left(  x\right)  :\left\Vert x\right\Vert \geq1\right\}
>0=Q\left(  0\right)  ,$

\item $Q\left(  x\right)  <4\rho\Rightarrow\left\Vert x\right\Vert <8\rho
\ $when $\rho\geq1$; otherwise $Q\left(  x\right)  <4\rho\Rightarrow\left\Vert
x\right\Vert <\delta\left(  \rho\right)  \equiv\left(  \left(  1+4\rho\right)
^{2n}-1\right)  ^{1/2n},$ this latter implication simply using $\left(
2.1\right)  $ and the definition of $Q$.

\item there exists $\delta>0$ such that $Q$ extends to a Lipschitz,
holomorphic map $\widetilde{Q}$ on the uniform strip $X\subset W_{\delta
}\subset\widetilde{X}$ given by,%

\[
W_{\delta}=\left\{  x+z:x\in X,\ z\in\widetilde{X},\ \left\Vert z\right\Vert
_{\widetilde{X}}<\delta\right\}  .
\]

\end{enumerate}

\medskip

\noindent We use the notion of a $Q$\textit{-body}, which for $\rho>0$ is
defined by%

\[
D_{Q}\left(  x,\rho\right)  =\left\{  y\in X:Q\left(  y-x\right)
<\rho\right\}  .
\]

\noindent Let $\left\Vert \cdot\right\Vert _{c_{0}}$ be an equivalent analytic
norm on $c_{0},$ with $\left\Vert x\right\Vert _{\infty}\leq\left\Vert
x\right\Vert _{c_{0}}\leq A_{1}\left\Vert x\right\Vert _{\infty}$ for all
$x\in c_{0},$ and some $A_{1}>1$ (see e.g., \cite{FPWZ}, and also \cite{AFK1},
\cite{FK} where it is referred to as the Preiss norm). It may be worth
pointing out that the Preiss norm $\left\Vert \cdot\right\Vert _{c_{0}}$ is
obtained as the restriction of a holomorphic map $\widetilde{\lambda}$ defined
on a neighbourhood of $c_{0}\backslash\left\{  0\right\}  $ in $\widetilde
{c}_{0}.$

\medskip

\noindent For the remainder of the proof, we fix a dense sequence $\left\{
x_{n}\right\}  _{n=1}^{\infty}$ in $X.$

\medskip

\begin{lemma}
\label{lemma phi}Let $\widetilde{V}=W_{\delta}$ be an open strip around $X$ in
$\widetilde{X}$ in which the function $\widetilde{Q}$ given above is defined.
Given $\varepsilon\in(0,1),$ $r>0,$ and a covering $\left\{  D_{Q}%
(x_{n},r)\right\}  _{n=1}^{\infty}$ of $X$, there exists a sequence of
holomorphic functions $\widetilde{\varphi}_{n}=\widetilde{\varphi
}_{n,r,\varepsilon}:\widetilde{V}\rightarrow\mathbb{C}$, whose restrictions to
$X$ we denote by $\varphi_{n}=\varphi_{n,r,\varepsilon}$, with the following properties:

\begin{description}
\item[1] The collection $\{\varphi_{n,r,\varepsilon}:X\rightarrow
\lbrack0,2]\,|\,n\in\mathbb{N}\}$ is equi-Lipschitz on $X$, with Lipschitz
constant of the form $L_{\varphi}=L_{1}\text{Lip}(Q)/r>1$ (where $L_{1}\geq1$
is independent of $\varepsilon$ and $n$),

\item[2] $0\leq\varphi_{n,r,\varepsilon}(x)\leq1+\varepsilon$ for all $x\in X$.

\item[3] For each $x\in X$ there exists $m=m_{x,r}\in\mathbb{N}$ (independent
of $\varepsilon)$ with $\varphi_{m,r,\varepsilon}(x)>1/2$.

\item[4] $0\leq\varphi_{n,r,\varepsilon}(x)\leq\varepsilon$ for all $x\in
X\setminus D_{Q}(x_{n},5r)$.

\item[4$^{\prime}$] $\left\Vert \varphi_{n,r,\varepsilon}^{\prime
}(x)\right\Vert \leq\varepsilon$ for all $x\in X\setminus D_{Q}(x_{n},5r)$.

\item[5] For each $x\in X$ there exist $\delta_{x,r}>0$, $a_{x,r}>0,\ $and
$n_{x,r}\in\mathbb{N}$ (all independent of $\varepsilon$) such that
\[
|\widetilde{\varphi}_{n,r,\varepsilon}(x+z)|<\frac{1}{M_{n}n!a_{x,r}^{n}%
}\,\,\text{ for }\,n>n_{x,r},\ \,z\in\widetilde{X}\text{ with }\Vert
z\Vert_{\widetilde{X}}<\delta_{x,r},
\]

where $M_{n}=e^{2C^{2}\kappa}\left(  1+\left\Vert x_{n}\right\Vert \right)  ,$
and the $\kappa=\kappa\left(  r\right)  >1$ and $C>1$ are constants that will
be specified in the proof of Theorem 1.

\item[6] For each $x\in X,$ there exists $\delta_{x,r}>0$ (independent of
$\varepsilon$) and $n_{x,\varepsilon,r}\in\mathbb{N}$ such that for $\Vert
z\Vert_{\widetilde{X}}<\delta_{x,r}$ and $n>n_{x,\varepsilon,r}$ we have
$|\widetilde{\varphi}_{n,r,\varepsilon}(x+z)|<\varepsilon$.

\item[7] For each $x\in X,$ there exists $\delta_{x,\varepsilon,r}$ such that
\[
|\widetilde{\varphi}_{n,r,\varepsilon}(x+z)|\leq1+2\varepsilon\,\,\text{ for
}n\in\mathbb{N},\ \text{and }z\in\widetilde{X}\text{ with }\Vert
z\Vert_{\widetilde{X}}\leq\delta_{x,\varepsilon,r}.
\]

\end{description}
\end{lemma}

\medskip

\noindent\textbf{Proof.\ }We largely follow the proof of \cite[Lemma 3]{AFK1},
with the few noted changes. As the proof in \cite{AFK1} is rather long and
technical, we here indicate only the key constructions, referring the reader
to the above cited paper for full details. Note that because $r$ is fixed
throughout, for ease of notation, we shall often suppress dependences on $r.$
Define subsets $A_{1,r}=\{y_{1}\in\mathbb{R}:-1\leq y_{1}\leq4r\}$, and, for
$n\geq2$,
\begin{align*}
A_{n,r}  &  =\{y=\{y_{j}\}_{j=1}^{n}\in\ell_{\infty}^{n}:-1-r\leq y_{n}%
\leq4r,\,2r\leq y_{j}\\
&  \leq M_{n,r}+2r\text{ for }1\leq j\leq n-1\},
\end{align*}%
\begin{align*}
A_{n,r}^{\prime}  &  =\{y=\{y_{j}\}_{j=1}^{n}\in\ell_{\infty}^{n}:-1\leq
y_{n}\leq3r,\,3r\leq y_{j}\\
&  \leq M_{n,r}+r\text{ for }1\leq j\leq n-1\},
\end{align*}%
\[
\text{where }M_{n,r}=\sup\left\{  Q\left(  x-x_{j}\right)  :x\in D_{Q}%
(x_{n},4r),\ 1\leq j\leq n\right\}  .
\]

\medskip

\noindent Let $\mu\in C^{\infty}\left(  \mathbb{R},\left[  0,1+\varepsilon
\right]  \right)  $ be Lipschitz such that $\mu\left(  t\right)  =0$ iff
$t\geq1,$ and $\mu\left(  t\right)  =1+\varepsilon$ iff $t\leq1/2.$ Let
$b^{n}\in C^{\infty}\left(  \mathbb{R},\left[  0,1\right]  \right)  $ be
Lipschitz such that $b^{n}\left(  t\right)  =1$ iff $t\notin\left(
2r,M_{n,r}+2r\right)  ,$ and $b^{n}\left(  t\right)  =0$ iff $t\in\left[
3r,M_{n,r}+r\right]  .$ Let $\widehat{b}\in C^{\infty}\left(  \mathbb{R}%
,\left[  0,1\right]  \right)  $ be Lipschitz such that $\widehat{b}\left(
t\right)  =1$ iff $t\notin\left(  -1-r,4r\right)  ,$ and $\widehat{b}\left(
t\right)  =0$ iff $t\in\left[  -1,3r\right]  .$ Now define a Lipschitz,
$C^{\infty}\ $smooth map $b_{n}:c_{00}\subset c_{0}\rightarrow\left[
0,1\right]  $ by $b_{n}\left(  y_{1},...,y_{n}\right)  =\mu\left(  \left\Vert
\left(  b^{n}\left(  y_{1}\right)  ,...,b^{n}\left(  y_{n-1}\right)
,\widehat{b}\left(  y_{n}\right)  \right)  \right\Vert _{c_{0}}\right)  .$
Then support$\left(  b_{n}\right)  =\overline{A}_{n},$ and $b_{n}%
=1+\varepsilon$ on $A_{n}^{\prime}.$ Moreover, $b_{n}$ is Lipschitz with
constant of the form $L_{1}/r,$ where $L_{1}\geq1$ is independent of $n.$

\medskip

\noindent Now one defines, on $\mathbb{R}^{n},$ the map%
\begin{align*}
h_{n}\left(  x\right)   &  =\frac{1}{T_{n}}\int_{\mathbb{R}^{n}}%
b_{n}(y)e^{-k\sum_{j=1}^{n}2^{-j}(x_{j}-y_{j})^{2}}dy,\\
& \\
T_{n}  &  =\int_{\mathbb{R}^{n}}e^{-k\sum_{j=1}^{n}2^{-j}y_{j}{}^{2}}dy.
\end{align*}

\medskip

\noindent Because $b_{n}=b_{n,\varepsilon}$ has compact support, is bounded,
Lipschitz, and $C^{1},$ one can choose $k=k_{n,\varepsilon}>0$ sufficiently
large that%

\begin{equation}
|b_{n}(x)-h_{n}(x)|\leq\varepsilon/2\,\,\text{ for all }\,\,x\in\mathbb{R}%
^{n},
\end{equation}

\noindent and%
\begin{equation}
|b_{n}^{\prime}(x)-h_{n}^{\prime}(x)|\leq\varepsilon/2\,\,\text{ for all
}\,\,x\in\mathbb{R}^{n}.
\end{equation}

\medskip

\noindent Next one defines (real) analytic maps $\varphi_{n}:X\rightarrow
\mathbb{R}$ by,
\[
\varphi_{n}\left(  x\right)  =h_{n}\left(  Q\left(  x-x_{1}\right)
,...,Q\left(  x-x_{n}\right)  \right)  =\frac{1}{T_{n}}\int_{\mathbb{R}^{n}%
}b_{n}\left(  y\right)  e^{-k_{n}\sum_{j=1}^{n}2^{-j}\left(  Q\left(
x-x_{j}\right)  -y_{j}\right)  ^{2}}.
\]

\noindent It is more or less standard to show that Lip$\left(  \varphi
_{n}\right)  \leq\frac{L_{1}}{r}$ Lip$\left(  Q\right)  .$ We can extend the
maps $\varphi_{n,r,\varepsilon}$ to complex valued maps defined on $W_{_{Q}}$
(see above) calling them $\widetilde{\varphi}_{n}$. Namely (where $x\in X$,
$z\in\widetilde{X}$),
\[
\widetilde{\varphi}_{n}\left(  x+z\right)  =\frac{1}{T_{n}}\int_{\mathbb{R}%
^{n}}b_{n}\left(  y\right)  e^{-k_{n}\sum_{j=1}^{n}2^{-j}\left(  \widetilde
{Q}\left(  x-x_{j}+z\right)  -y_{j}\right)  ^{2}}dy
\]
\noindent Note that the $\widetilde{\varphi}_{n}$ are well defined (as the
$b_{n}$ have compact supports) and are holomorphic where $\widetilde{Q}$ is
(namely on $\widetilde{W}_{\delta}$).

\medskip

\noindent To see $\left(  4\right)  $ and $\left(  4^{\prime}\right)  $, note
that if $Q\left(  x-x_{n}\right)  \geq5r,$ then there is a neighbourhood $N$
of $x$ for which $y\in N$ implies that the point,
\[
\widehat{x}=\left(  Q\left(  y-x_{1}\right)  ,...,Q\left(  y-x_{n}\right)
\right)  \in\mathbb{R}^{n}\backslash A_{n},
\]
from which we have $b_{n}\left(  \widehat{x}\right)  =0$ and $b_{n}^{\prime
}\left(  \widehat{x}\right)  =0.$ Hence, by $\left(  2.2\right)  $ and
$\left(  2.3\right)  ,$ we have, $\left\vert \varphi_{n}\left(  x\right)
\right\vert <\varepsilon/2$ and $\left\Vert \varphi_{n}^{\prime}\left(
x\right)  \right\Vert <\varepsilon/2.$

\medskip

\noindent The remaining parts are handled as in \cite{AFK1}, noting that for
$\left(  5\right)  $ we choose $\kappa_{n}$ larger if need be to ensure the
stated estimate involving the $M_{n}.$ $\square$

\bigskip

\noindent We return now to the proof of the theorem. Let $\varepsilon>0$ be
given and choose $\varepsilon^{\prime}$ satisfying
\[
0<\varepsilon^{\prime}<\min\{\frac{1}{8},1/(132C_{0}A_{1}^{2}L_{1}%
Lip(Q)),1/(10A_{1}r)\},
\]
where $L_{1}$ is as in part (1) of the preceding lemma, where \ is defined
immediately below, and where $C_{0}$ is a constant, only depending on $X$,
which will be fixed later on (see page 9 below). Because $f$ is bounded, we
may suppose that $1\leq f\leq2.$ As $f^{\prime}$ is uniformly continuous on
$X$, we can find a fixed $\rho>0$ so that for all $n,$ $x\in B_{\rho}\left(
x_{n}\right)  $ implies $\left\Vert f^{\prime}\left(  x_{n}\right)
-f^{\prime}\left(  x\right)  \right\Vert <\varepsilon^{\prime}.$ Now,
considering property $\left(  4\right)  $ of $Q,$ and noting that
$\delta\left(  r\right)  \rightarrow0^{+}$ as $r\rightarrow0^{+},$ we can
choose $r\in\left(  0,1\right)  $ (independent of $n$) so that $D_{Q}\left(
x_{n},5r\right)  \subset B_{\rho}\left(  x_{n}\right)  $ for all $n.$ It will
be convenient to write $D_{n}\equiv D_{Q}\left(  x_{n},5r\right)  .$ This $r$
shall be fixed for the remainder of the proof.

\medskip

\subsection{The functions $\nu_{n}$}

Next let $\overline{\nu}\in C^{\infty}\left(  \mathbb{R},\left[  0,1\right]
\right)  $ be Lipschitz such that $\overline{\nu}\left(  t\right)  =1$ iff
$\left\vert t\right\vert \leq5r,$ and $\overline{\nu}\left(  t\right)  =0$ iff
$\left\vert t\right\vert \geq\frac{11}{2}r.$ Put $L=$ Lip$\left(  f\right)  $.
Fix a sequence of functions $\left\{  \varphi_{n,r,\varepsilon_{1}}\right\}
_{n=1}^{\infty}$ with respect to the covering $\left\{  D_{Q}\left(
x_{n},r\right)  \right\}  _{n=1}^{\infty}$ of $X$ as given by Lemma
\ref{lemma phi}, where $r$ is fixed as above and the $\varepsilon$ of the
Lemma is chosen to be
\[
\varepsilon_{1}:=\min\left\{  \varepsilon^{\prime}r/3C_{0}L\text{Lip}\left(
\overline{\nu}\right)  ,\varepsilon^{\prime}r/25L\text{Lip}\left(
\overline{\nu}\right)  \right\}  ,
\]

\noindent We write $\varphi_{n,r,\varepsilon_{1}}$ as $\varphi_{n}$ for ease
of notation, and, as in Lemma \ref{lemma phi} $\left(  1\right)  ,$
$L_{\varphi}=$ Lip$\left(  \varphi_{n}\right)  =L_{1}$Lip$\left(  Q\right)
/r\geq1$, which we recall is independent of $n.$ Often we will subsume
dependence on $\varepsilon_{1}$ as dependence on $\varepsilon^{\prime}$ and
$L.$

\medskip

\noindent Put $\Delta\left(  t\right)  =\left(  \left(  \left\vert
t\right\vert +1\right)  ^{2n}-1\right)  ^{1/2n}\geq0.$ Now via convolution
integrals between $\overline{\nu}$ and Gaussian kernels, we can find
Lipschitz, analytic functions $\nu,$ with Lip$\left(  \nu\right)  =$
Lip$\left(  \overline{\nu}\right)  ,$ and which $C^{1}$-fine approximate
$\overline{\nu}$ in the following sense,

\medskip%

\begin{align}
\left\vert \nu\left(  t\right)  -\overline{\nu}\left(  t\right)  \right\vert
&  <\frac{\varepsilon^{\prime}r/2LL_{\varphi}}{1+\Delta\left(  t\right)
},\nonumber\\
& \\
\left\vert \nu^{\prime}\left(  t\right)  -\overline{\nu}^{\prime}\left(
t\right)  \right\vert  &  <\frac{\varepsilon^{\prime}r/2LL_{\varphi}}%
{1+\Delta\left(  t\right)  }.\nonumber
\end{align}

\medskip

\noindent Indeed, we can take $\nu$ to be of the form,

\medskip%

\begin{align*}
\nu\left(  t\right)   &  =\frac{1}{a}\int_{\mathbb{R}}\overline{\nu}\left(
s\right)  e^{-\kappa\left(  t-s\right)  ^{2}}ds,\\
& \\
a  &  =\int_{\mathbb{R}}e^{-\kappa s^{2}}ds,
\end{align*}

\noindent where $\kappa>1$ is chosen sufficiently large and is independent of
$t$ (although it does depend on $\max\left\{  \Delta\left(  t\right)
:t\in\ \text{supp}\left(  \overline{\nu}\right)  \right\}  <\infty$). This is
possible because $\overline{\nu}$ is $C^{\infty}$ with compact support, and
the function $t\rightarrow\frac{\varepsilon^{\prime}/2L}{1+\Delta\left(
t\right)  }$ is strictly positive, continuous and decreases slowly enough with
respect to $e^{-\kappa t^{2}}$ (namely, $\lim_{t\rightarrow\infty}%
\Delta(t)/e^{\kappa t^{2}}=0$). Moreover, since $\overline{\nu}$ has compact
support, $\nu$ has a holomorphic extension,
\[
\widetilde{\nu}\left(  z\right)  =\frac{1}{a}\int_{\mathbb{R}}\overline{\nu
}\left(  s\right)  e^{-\kappa\left(  z-s\right)  ^{2}}ds,
\]
to $\mathbb{C}.$ Next observe that for $t,s\in\mathbb{R}$ and $z\in\mathbb{C}$
with $\left\vert z\right\vert \leq\eta,$ we have,%

\begin{align*}
\operatorname{Re}\left(  t+z-s\right)  ^{2}  &  =\left(  t-s\right)
^{2}+2\left(  t-s\right)  \operatorname{Re}z+\operatorname{Re}\left(
z^{2}\right) \\
& \\
&  =\left(  t-s+\operatorname{Re}z\right)  ^{2}-\left(  \operatorname{Re}%
z\right)  ^{2}+\operatorname{Re}\left(  z^{2}\right) \\
& \\
&  \geq\left(  t-s+\operatorname{Re}z\right)  ^{2}-2\eta^{2}.
\end{align*}

\noindent Therefore when $\left\vert z\right\vert <\eta$ we get,
\begin{align}
\left\vert \widetilde{\nu}\left(  t+z\right)  \right\vert  &  =\frac{1}%
{a}\left\vert \int_{\mathbb{R}}\overline{\nu}\left(  s\right)  e^{-\kappa
\left(  t+z-s\right)  ^{2}}ds\right\vert \nonumber\\
& \nonumber\\
&  \leq\frac{1}{a}\int_{\mathbb{R}} e^{-\kappa\operatorname{Re}\left(
t+z-s\right)  ^{2}}ds\nonumber\\
& \nonumber\\
&  \leq\frac{1}{a}\int_{\mathbb{R}} e^{-\kappa\left(  t-s+\operatorname{Re}%
z\right)  ^{2}-2\eta^{2}}ds\\
& \nonumber\\
&  =\frac{e^{2\kappa\eta^{2}}}{a}\int_{\mathbb{R}} e^{-\kappa\left(
t+\operatorname{Re}z-s\right)  ^{2}}ds\nonumber\\
& \nonumber\\
&  =e^{2\kappa\eta^{2}},\nonumber
\end{align}

\medskip

\noindent where we have used a variable change to obtain the last line. Now
define Lipschitz, analytic functions $\nu_{n}:X\rightarrow\mathbb{R}$ by,
\[
\nu_{n}\left(  x\right)  =\nu\left(  Q\left(  x-x_{n}\right)  \right)  .
\]
Clearly $\nu_{n}$ has the holomorphic extension $\widetilde{\nu}_{n}\left(
z\right)  =\widetilde{\nu}\left(  \widetilde{Q}\left(  z-x_{n}\right)
\right)  . $ It will be convenient to put $\overline{\nu}_{n}\left(  x\right)
=\overline{\nu}\left(  Q\left(  x-x_{n}\right)  \right)  .$ Observe that,
writing $\widehat{D}_{n}=D_{Q}\left(  x_{n},6r\right)  ,$%

\begin{equation}
\left\vert \nu_{n}\left(  x\right)  \right\vert <\frac{\varepsilon^{\prime
}r/2LL_{\varphi}}{1+\Delta\left(  Q\left(  x-x_{n}\right)  \right)
},\ \ \text{for\ }x\notin\widehat{D}_{n},
\end{equation}

and%

\begin{equation}
\left\vert \nu_{n}^{\prime}\left(  x\right)  \right\vert <\frac{\text{Lip}%
\left(  Q\right)  \varepsilon^{\prime}r/2LL_{\varphi}}{1+\Delta\left(
Q\left(  x-x_{n}\right)  \right)  },\ \ \text{for\ }x\notin\widehat{D}_{n}.
\end{equation}

\noindent Note that
\begin{align}
\frac{\varepsilon^{\prime}r/2LL_{\varphi}}{1+\Delta\left(  Q\left(
x-x_{n}\right)  \right)  }  &  =\frac{\varepsilon^{\prime}r/2LL_{\varphi}%
}{1+q\left(  x-x_{n}\right)  ^{1/2n}}\nonumber\\
& \\
&  \leq\frac{\varepsilon^{\prime}r/2LL_{\varphi}}{1+\left\Vert x-x_{n}%
\right\Vert }.\nonumber
\end{align}

\medskip

\noindent Now we estimate $\left\vert \widetilde{\nu}_{n}\left(  x+z\right)
\right\vert =\left\vert \widetilde{\nu}\left(  \widetilde{Q}\left(
x-x_{n}+z\right)  \right)  \right\vert ,$ for $\left\Vert z\right\Vert
_{\widetilde{X}}<\eta.$ From \cite[Lemma 2]{AFK1}, we can write%

\[
\widetilde{Q}\left(  x-x_{n}+z\right)  =Q\left(  x-x_{n}\right)  +Z_{n},
\]

\noindent where $Z_{n}\in\mathbb{C}$ with $\left\vert Z_{n}\right\vert \leq
C\left\Vert z\right\Vert _{\widetilde{X}},$ for some constant $C>1.$ Then from
the calculation $\left(  2.5\right)  $ we get, for $\left\Vert z\right\Vert
_{\widetilde{X}}<\eta,$%

\begin{align}
\left\vert \widetilde{\nu}_{n}\left(  x+z\right)  \right\vert  &  =\left\vert
\widetilde{\nu}\left(  \widetilde{Q}\left(  x-x_{n}+z\right)  \right)
\right\vert \nonumber\\
& \nonumber\\
&  =\left\vert \widetilde{\nu}\left(  Q\left(  x-x_{n}\right)  +Z_{n}\right)
\right\vert \\
& \nonumber\\
&  \leq e^{2C^{2}\kappa\eta^{2}}.\nonumber
\end{align}

\medskip

\noindent It is also worthwhile to note that $\nu\left(  t\right)
<1+\varepsilon^{\prime}$ for all $t.$

\medskip

\noindent Let $T_{n}\left(  x\right)  =f^{\prime}\left(  x_{n}\right)  \left(
x-x_{n}\right)  +f\left(  x_{n}\right)  $ be the first order Taylor polynomial
of $f$ at $x_{n}.$ Note that $\left\Vert T_{n}^{\prime}\left(  x\right)
\right\Vert =\left\Vert f^{\prime}\left(  x_{n}\right)  \right\Vert \leq L.$
Observe that $T_{n}-f$ is Lipschitz on $B_{\rho}\left(  x_{n}\right)  ,$ with
Lip$\left(  T_{n}-f\right)  \leq\left\Vert \left(  T_{n}-f\right)  ^{\prime
}\right\Vert =\left\Vert f^{\prime}\left(  x_{n}\right)  -f^{\prime}\left(
x\right)  \right\Vert \leq\varepsilon^{\prime}$ on $B_{\rho}\left(
x_{n}\right)  .$ It follows that $T_{n}-f$ is Lipschitz on $D_{n}\subset
B_{\rho}\left(  x_{n}\right)  $ with constant not greater than $\varepsilon
^{\prime}.$ Denote by $\overline{T_{n}-f}$ a bounded and Lipschitz extension
of $\left(  T_{n}-f\right)  \mid_{D_{n}}$ to all of $X,$ having the same bound
and Lipschitz constant. For example, one can take, temporarily writing
$h=\left(  T_{n}-f\right)  \mid_{D_{n}},$%
\[
\left(  \overline{T_{n}-f}\right)  \left(  x\right)  =\max\{-\Vert
h\Vert_{\infty},\min\{\Vert h\Vert_{\infty},\,\inf_{y\in D_{n}}%
\{h(y)+\text{Lip}(h)\Vert x-y\Vert\}\,\}\,\}.
\]

\medskip

\noindent Write $\epsilon_{n}\left(  x\right)  =\left(  \overline{T_{n}%
-f}\right)  \left(  x\right)  .$ We now apply \cite[Proposition 3]{AFK1} to
$\epsilon_{n}\left(  x\right)  ,$ along with the standard `scaling argument'
that appears at the very end of the proof of \cite[Theorem 1]{AFK1}, to obtain
the following: there exists a constant $C_{0}>1,$ depending only on $X,$ a
neighbourhood $X\subset\widetilde{W}\subset\widetilde{X},$ where
$\widetilde{W}=\widetilde{W}_{\varepsilon^{\prime},r}$ depends only on
$\varepsilon^{\prime}$ and $r$ (the dependence on $L_{\varphi}$ written as a
dependence on $r$), and an analytic map $\delta_{n}:X\rightarrow\mathbb{R}$
such that

\begin{enumerate}
\item $\left\vert \epsilon_{n}\left(  x\right)  -\delta_{n}\left(  x\right)
\right\vert <\varepsilon^{\prime}r/L_{\varphi}$ for all $x\in X,$

\item Lip$\left(  \delta_{n}\right)  \leq C_{0}$ Lip$\left(  \epsilon
_{n}\right)  \leq C_{0}\varepsilon^{\prime},$

\item the map $\delta_{n}$ extends to a holomorphic map $\widetilde{\delta
}_{n}$ on $\widetilde{W}$ (where in particular, $\widetilde{W}$ is independent
of $n$)$,$

\item $|\widetilde{\delta}_{n}(x+iy)-\delta_{n}(x)|\leq M_{\Delta}$ for all
$x+iy\in\widetilde{W},$ where $M_{\Delta}$ depends on $\varepsilon^{\prime}$
and is independent of $n.$
\end{enumerate}

\medskip

\noindent Now we define analytic functions on $X$ by,%

\[
\psi_{n}\left(  x\right)  =\left(  T_{n}\left(  x\right)  \nu_{n}\left(
x\right)  -\delta_{n}\left(  x\right)  \right)  \varphi_{n}\left(  x\right)
.
\]

\medskip

\noindent Observe that from property $\left(  3\right)  $ of $\delta_{n}$ and
Lemma \ref{lemma phi}, $\psi_{n}$ extends to a holomorphic map $\widetilde
{\psi}_{n}\left(  z\right)  =\left(  \widetilde{T}_{n}\left(  z\right)
\widetilde{\nu}_{n}\left(  z\right)  -\widetilde{\delta}_{n}\left(  z\right)
\right)  \widetilde{\varphi}_{n}\left(  z\right)  ,$ where
\[
\widetilde{T}_{n}\left(  z\right)  =\widetilde{T}_{n}\left(  x+iy\right)
=\widetilde{f^{\prime}\left(  x_{n}\right)  }\left(  x+iy-x_{n}\right)
+f\left(  x_{n}\right)
\]
($\widetilde{f^{\prime}\left(  x_{n}\right)  }$ being the canonical extension
of $f^{\prime}\left(  x_{n}\right)  $ to all of $\widetilde{X}$), on a
neighbourhood $X\subset\widetilde{W}\subset\widetilde{X},$ where
$\widetilde{W}$ is independent of $n.$

\medskip

\subsection{The approximating function $g$}

Let us define the function $g:X\rightarrow\mathbb{R}$ by,%
\[
g\left(  x\right)  =\frac{\left\Vert \left\{  \psi_{n}\left(  x\right)
\right\}  _{n=1}^{\infty}\right\Vert _{c_{0}}}{\left\Vert \left\{  \varphi
_{n}\left(  x\right)  \right\}  _{n=1}^{\infty}\right\Vert _{c_{0}}}%
\]

\medskip

\noindent We next show that $g$ is analytic. Since the norm $\|\cdot\|_{c_{0}%
}$ is real analytic on $c_{0}\setminus\{0\}$, it is sufficient to check that
the mappings $x\mapsto\{\varphi_{n}(x)\}_{n=1}^{\infty}$ and $x\mapsto
\{\psi_{n}(x)\}_{n=1}^{\infty}$ are real analytic from $X$ into $c_{0}$ and do
not take the value $0\in c_{0}$. Using Lemma \ref{lemma phi} it is easy to
show that the function $x\mapsto\{\varphi_{n}(x)\}_{n=1}^{\infty}$ has such
properties (see \cite[Lemma 4]{AFK1}).

As for the function $x\mapsto\{\psi_{n}(x)\}_{n=1}^{\infty}$, let us first
show that it does not take the value $0$. In fact we show that for each $x\in
X$ there exists an $n$ so that the number $\left(  T_{n}\left(  x\right)
\nu_{n}\left(  x\right)  -\delta_{n}\left(  x\right)  \right)  \varphi
_{n}\left(  x\right)  $ is bounded above $1/4$. Indeed, for each $x\in X,$
there is a minimal $n=n_{x}$ with $x\in D_{Q}\left(  x_{n_{x}},3r\right)  ,$
and via the proof of \cite[Lemma 3 $\left(  3\right)  $]{AFK1}, for such
$n_{x}$ we have $\varphi_{n_{x}}\left(  x\right)  \geq1/2.$ Note also that
$D_{Q}\left(  x_{n_{x}},3r\right)  \subset D_{Q}\left(  x_{n_{x}},5r\right)
=D_{n_{x}},$ and $\epsilon_{n_{x}}\left(  x\right)  =T_{n_{x}}\left(
x\right)  -f\left(  x\right)  $ on $D_{n_{x}}.$ So, from this and property
$\left(  1\right)  $ of $\delta_{n},$ we have $\left\vert T_{n_{x}}\left(
x\right)  \overline{\nu}_{n_{x}}\left(  x\right)  -f\left(  x\right)
-\delta_{n_{x}}\left(  x\right)  \right\vert =\left\vert \epsilon_{n_{x}%
}\left(  x\right)  -\delta_{n_{x}}\left(  x\right)  \right\vert \leq
\varepsilon^{\prime}.$ Now to replace $\overline{\nu}_{n_{x}}$ with
$\nu_{n_{x}},$ we observe by $\left(  2.4\right)  \ $and $\left(  2.8\right)
,$%

\begin{align}
\left\vert T_{n}\left(  x\right)  \nu_{n}\left(  x\right)  -T_{n}\left(
x\right)  \overline{\nu}_{n}\left(  x\right)  \right\vert  &  =\left\vert
T_{n}\left(  x\right)  \right\vert \left\vert \nu_{n}\left(  x\right)
-\overline{\nu}_{n}\left(  x\right)  \right\vert \nonumber\\
& \nonumber\\
&  \leq\left(  L\left\Vert x-x_{n}\right\Vert +\left\vert f\left(
x_{n}\right)  \right\vert \right)  \ \frac{\varepsilon^{\prime}r/2LL_{\varphi
}}{1+\left\Vert x-x_{n}\right\Vert }\nonumber\\
& \\
&  \leq\left(  L\left\Vert x-x_{n}\right\Vert +2\right)  \ \frac
{\varepsilon^{\prime}r/2LL_{\varphi}}{1+\left\Vert x-x_{n}\right\Vert
}\nonumber\\
& \nonumber\\
&  \leq\varepsilon^{\prime}r/2L_{\varphi}+\varepsilon^{\prime}r/LL_{\varphi
}\nonumber\\
& \nonumber\\
&  \leq3\varepsilon^{\prime}r/L_{\varphi}\leq3\varepsilon^{\prime}.\nonumber
\end{align}

\noindent Therefore, these estimates give, $\left\vert T_{n_{x}}\left(
x\right)  \nu_{n_{x}}\left(  x\right)  -f\left(  x\right)  -\delta_{n_{x}%
}\left(  x\right)  \right\vert \leq4\varepsilon^{\prime},$ and because
$f\geq1,$ we have our desired bound
\begin{align*}
\left\vert T_{n_{x}}\left(  x\right)  \nu_{n_{x}}\left(  x\right)
-\delta_{n_{x}}\left(  x\right)  \right\vert \varphi_{n_{x}}\left(  x\right)
&  \geq\left\vert T_{n_{x}}\left(  x\right)  \nu_{n_{x}}\left(  x\right)
-\delta_{n_{x}}\left(  x\right)  \right\vert \left(  1/2\right) \\
& \\
&  \geq\left(  f\left(  x\right)  -4\varepsilon^{\prime}\right)  \left(
1/2\right)  >1/4.
\end{align*}
We remark that it follows from this that for any $x,$
\begin{equation}
\left\Vert \left\{  \psi_{n}\left(  x\right)  \right\}  _{n}\right\Vert
_{c_{0}}\geq\left\Vert \left\{  \psi_{n}\left(  x\right)  \right\}
_{n}\right\Vert _{\infty}=\left\Vert \left\{  \left(  T_{n}\left(  x\right)
-\delta_{n}\left(  x\right)  \right)  \varphi_{n}\left(  x\right)  \right\}
_{n}\right\Vert _{\infty}\geq1/4.
\end{equation}

\medskip

\noindent Next, to show that the function $x\mapsto\{\psi_{n}(x)\}_{n=1}%
^{\infty}$ is real analytic from $X$ into $c_{0}$, we shall require that for
each $x,$ there exists $n_{x}$ and $\delta_{x}\in\left(  0,1\right)  $ so that
for $n\geq n_{x}$ and $\left\Vert z\right\Vert _{\widetilde{X}}<\delta_{x},$
we have
\begin{equation}
\frac{\left\vert \widetilde{T}_{n}\left(  x+z\right)  \widetilde{\nu}%
_{n}\left(  x+z\right)  -\widetilde{\delta}_{n}\left(  x+z\right)  \right\vert
}{M_{n}}\leq M_{x},
\end{equation}
where $M_{x}$ depends on $x,$ but is independent of $n.$

\medskip

\noindent Recalling that $\widetilde{T}_{n}\left(  w\right)  =$ $\widetilde
{f^{\prime}\left(  x_{n}\right)  }\left(  x-x_{n}+w\right)  +f\left(
x_{n}\right)  ,$ and using $\left(  2.9\right)  $ and property $\left(
1\right)  \ $and $\left(  4\right)  $ of $\delta_{n},$ where we may suppose
that $x+z\in\widetilde{W}$ when $\left\Vert z\right\Vert _{\widetilde{X}%
}<\delta_{x}<1,$ we obtain,
\begin{align*}
&  \left\vert \widetilde{T}_{n}\left(  x+z\right)  \widetilde{\nu}_{n}\left(
x+z\right)  -\widetilde{\delta}_{n}\left(  x+z\right)  \right\vert \\
& \\
&  \leq\left\vert \widetilde{T}_{n}\left(  x+z\right)  \right\vert \left\vert
\widetilde{\nu}_{n}\left(  x+z\right)  \right\vert +\left\vert \widetilde
{\delta}_{n}\left(  x+z\right)  -\delta_{n}\left(  x\right)  \right\vert
+\delta_{n}\left(  x\right) \\
& \\
&  \leq\left\vert \widetilde{f^{\prime}\left(  x_{n}\right)  }\left(
x-x_{n}+z\right)  +f\left(  x_{n}\right)  \right\vert e^{2C^{2}\kappa
\delta_{x}^{2}}+M_{\Delta}+\varepsilon^{\prime}\\
& \\
&  <\left(  L\left(  \left\Vert x-x_{n}\right\Vert +\left\Vert z\right\Vert
_{\widetilde{X}}\right)  +\left\vert f\left(  x_{n}\right)  \right\vert
\right)  e^{2C^{2}\kappa\delta_{x}^{2}}+2M_{\Delta}\\
& \\
&  <\left(  L\left(  \left\Vert x-x_{n}\right\Vert +1\right)  +2\right)
e^{2C^{2}\kappa}+2M_{\Delta}\\
& \\
&  \leq\left(  3L\left(  \left\Vert x-x_{n}\right\Vert +1\right)  \right)
e^{2C^{2}\kappa}+2M_{\Delta}%
\end{align*}

\medskip

\noindent Now recalling that $M_{n}=e^{2C^{2}\kappa}\left(  1+\left\Vert
x_{n}\right\Vert \right)  ,$ we see that
\begin{align*}
&  \frac{3L\left(  \left\Vert x-x_{n}\right\Vert +1\right)  e^{2C^{2}\kappa}%
}{M_{n}}\\
& \\
&  \leq\frac{3L\left(  \left\Vert x\right\Vert +\left\Vert x_{n}\right\Vert
+1\right)  }{1+\left\Vert x_{n}\right\Vert }\\
& \\
&  \leq3L\left(  1+\left\Vert x\right\Vert \right)  .
\end{align*}
Putting $M_{x}=2M_{\Delta}+3L\left(  1+\left\Vert x\right\Vert \right)  ,$ we
have established $\left(  2.12\right)  .$

\medskip

\noindent Now, to show the analyticity of $\left\{  \psi_{n}\left(  x\right)
\right\}  _{n=1}^{\infty}$, we first note that property $(5)$ of Lemma
\ref{lemma phi} together with $(2.12)$ yield
\[
|\widetilde{\psi}_{n}(x+z)|=|\widetilde{T}_{n}(x+z)\widetilde{\nu}%
_{n}(x+z)-\widetilde{\delta}_{n}(x+z)|\, |\widetilde{\varphi}_{n}%
(x+z)|\leq\frac{M_{x}}{n! a_{x,r}^{n}}
\]
whenever $n\geq n_{x}$ and $\|z\|_{\widetilde{X}}<\delta_{x}$.

Because the numerical series $\sum_{n=1}^{\infty} M_{x}/n! a_{x,r}^{n}$ is
convergent, we then have that the series of functions $\sum_{n=1}^{\infty}|
\widetilde{\psi}_{n}\left(  x+z\right)  | $ is uniformly convergent on the
ball $B_{\widetilde{X}}(0, \delta_{x})$, which clearly implies that the
series
\[
\sum_{n=1}^{\infty}\widetilde{\psi}_{n}(z)e_{n} =\{\widetilde{\psi}%
_{n}(z)\}_{n=1}^{\infty}
\]
is uniformly convergent for $z\in B_{\widetilde{X}}(x, \delta_{x})$. Then it
is clear that $\{\widetilde{\psi}_{n}(z)\}_{n=1}^{\infty}$, being a series of
holomorphic mappings which converges uniformly on the ball $B_{\widetilde{X}%
}(x, \delta_{x})$, is a holomorphic mapping on this ball. Since $x\in X$ is
arbitrary, this shows that $x\mapsto\{\psi_{n}(x)\}_{n=1}^{\infty}$ is real analytic.

\medskip

\subsection{Showing that $g$ does the job}

Now we move on to our final estimates; $\left\vert g-f\right\vert $ and
$\left\Vert g^{\prime}-f^{\prime}\right\Vert .$ Fix $x\in X,$ and put
$\mathcal{N}=\mathcal{N}_{x}=\left\{  n:x\in D_{n}\right\}  .$ Now we have
(using $f\geq1>0),$ that
\begin{align*}
\left\vert g\left(  x\right)  -f\left(  x\right)  \right\vert  &  =\left\vert
\frac{\left\Vert \left\{  \psi_{n}\left(  x\right)  \right\}  _{n=1}^{\infty
}\right\Vert _{c_{0}}}{\left\Vert \left\{  \varphi_{n}\left(  x\right)
\right\}  _{n=1}^{\infty}\right\Vert _{c_{0}}}-f\left(  x\right)  \right\vert
\\
& \\
&  =\left\vert \frac{\left\Vert \left\{  \psi_{n}\left(  x\right)  \right\}
_{n=1}^{\infty}\right\Vert _{c_{0}}}{\left\Vert \left\{  \varphi_{n}\left(
x\right)  \right\}  _{n=1}^{\infty}\right\Vert _{c_{0}}}-\frac{\left\Vert
\left\{  f\left(  x\right)  \varphi_{n}\left(  x\right)  \right\}
_{n=1}^{\infty}\right\Vert _{c_{0}}}{\left\Vert \left\{  \varphi_{n}\left(
x\right)  \right\}  _{n=1}^{\infty}\right\Vert _{c_{0}}}\right\vert \\
& \\
&  =\frac{1}{\left\Vert \left\{  \varphi_{n}\left(  x\right)  \right\}
_{n=1}^{\infty}\right\Vert _{c_{0}}}\left\Vert \left\{  \left(  T_{n}\left(
x\right)  \nu_{n}\left(  x\right)  -f\left(  x\right)  -\delta_{n}(x)\right)
\varphi_{n}\left(  x\right)  \right\}  _{n=1}^{\infty}\right\Vert _{c_{0}}\\
& \\
&  \leq2\left\Vert \left\{  \left(  T_{n}\left(  x\right)  \nu_{n}\left(
x\right)  -f\left(  x\right)  -\delta_{n}(x)\right)  \varphi_{n}\left(
x\right)  \right\}  _{n=1}^{\infty}\right\Vert _{c_{0}},
\end{align*}

\medskip

\noindent the last line by Lemma \ref{lemma phi}\ $\left(  3\right)  .$ We
proceed in cases.

\medskip

\textbf{Case 1:\ \ }For $n\in\mathcal{N},$ we have $\epsilon_{n}\left(
x\right)  =T_{n}\left(  x\right)  -f\left(  x\right)  =T_{n}\left(  x\right)
\overline{\nu}_{n}\left(  x\right)  -f\left(  x\right)  ,$ and so by property
$\left(  1\right)  $ of $\delta_{n}$ and Lemma \ref{lemma phi} $\left(
7\right)  ,$ we obtain the estimate, $\left\vert T_{n}\left(  x\right)
\overline{\nu}_{n}\left(  x\right)  -f\left(  x\right)  -\delta_{n}%
(x)\right\vert \varphi_{n}\left(  x\right)  \leq\left(  \varepsilon^{\prime
}r/L_{\varphi}\right)  3.$ Then using $\left(  2.10\right)  ,$ we have
\begin{equation}
\left\vert T_{n}\left(  x\right)  \nu_{n}\left(  x\right)  -f\left(  x\right)
-\delta_{n}(x)\right\vert \varphi_{n}\left(  x\right)  \leq6r\varepsilon
^{\prime}/L_{\varphi}.
\end{equation}

\medskip

\textbf{Case 2:\ \ }For $n\notin\mathcal{N},$ recall $\varphi_{n}\left(
x\right)  \leq\varepsilon_{1}\leq\varepsilon^{\prime}r/25L.$

\medskip

\noindent Now, for $n$ such that $x\in\widehat{D}_{n},$ we have $Q\left(
x-x_{n}\right)  <6r,$ and so $\left\Vert x-x_{n}\right\Vert \leq\left(
q\left(  x-x_{n}\right)  \right)  ^{1/2n}<\left(  \left(  6r+1\right)
^{2n}-1\right)  ^{1/2n}<7,$ as $r<1.$ Hence, for such $n$ we have,%

\begin{equation}
\left\vert T_{n}\left(  x\right)  \nu_{n}\left(  x\right)  \right\vert
\leq\left(  L\left\Vert x-x_{n}\right\Vert +\left\vert f\left(  x_{n}\right)
\right\vert \right)  \nu_{n}\left(  x\right)  \leq\left(  7L+2\right)  \left(
1+\varepsilon^{\prime}\right)  \leq18L.
\end{equation}

\medskip

\noindent On the other hand, for $n$ such that $x\notin\widehat{D}_{n},$ by
$\left(  2.6\right)  $ we obtain,%

\begin{align*}
\left\vert T_{n}\left(  x\right)  \nu_{n}\left(  x\right)  \right\vert  &
\leq\left(  L\left\Vert x-x_{n}\right\Vert +\left\vert f\left(  x_{n}\right)
\right\vert \right)  \frac{\varepsilon^{\prime}r/2LL_{\varphi}}{1+\left\Vert
x-x_{n}\right\Vert }\\
& \\
&  \leq\left(  L\left\Vert x-x_{n}\right\Vert +2\right)  \frac{\varepsilon
^{\prime}/2L}{1+\left\Vert x-x_{n}\right\Vert }\\
& \\
&  \leq\varepsilon^{\prime}/2+\varepsilon^{\prime}/L\leq2\varepsilon^{\prime}.
\end{align*}

\noindent In any event, for all $n$ we have,
\begin{equation}
\left\vert T_{n}\left(  x\right)  \nu_{n}\left(  x\right)  \right\vert
\leq18L.
\end{equation}

\medskip

\noindent Therefore, for $n\notin\mathcal{N},$ using again property $\left(
1\right)  $ of $\delta_{n},$ we have,
\begin{align}
\left\vert T_{n}\left(  x\right)  \nu_{n}\left(  x\right)  -f\left(  x\right)
-\delta_{n}(x)\right\vert \varphi_{n}\left(  x\right)   &  \leq\left(
\left\vert T_{n}\left(  x\right)  \nu_{n}\left(  x\right)  \right\vert
+\left\vert f\left(  x\right)  \right\vert +\delta_{n}(x)\right)  \varphi
_{n}\left(  x\right) \nonumber\\
& \\
&  \leq\left(  18L+2+2\varepsilon^{\prime}\right)  \left(  \varepsilon
^{\prime}r/25L\right)  \leq\varepsilon^{\prime}r.\nonumber
\end{align}

\medskip

\noindent It follows that, $\left\vert g\left(  x\right)  -f\left(  x\right)
\right\vert \leq10A_{1}\varepsilon^{\prime}r<\varepsilon.$

\medskip

\noindent We now establish some derivative estimates. Fix $x$ and consider the
expression
\[
\left(  T_{n}\left(  x\right)  \nu_{n}\left(  x\right)  \right)  ^{\prime
}=T_{n}^{\prime}\left(  x\right)  \nu_{n}\left(  x\right)  +T_{n}\left(
x\right)  \nu_{n}^{\prime}\left(  x\right)  .
\]
From an estimate analogous to $\left(  2.15\right)  ,$ using $\left(
2.4\right)  $ and $\left(  2.7\right)  ,$ we have that for all $n$,
$\left\Vert T_{n}\left(  x\right)  \nu_{n}^{\prime}\left(  x\right)
\right\Vert \leq9L$Lip$\left(  Q\right)  $Lip$\left(  \nu\right)  .$ Also,
$\left\Vert T_{n}^{\prime}\left(  x\right)  \nu_{n}\left(  x\right)
\right\Vert \leq L\left(  1+\varepsilon^{\prime}\right)  \leq2L.$ Hence,
\[
\left\Vert \left(  T_{n}\left(  x\right)  \nu_{n}\left(  x\right)  \right)
^{\prime}\right\Vert \leq2L+9L\text{Lip}\left(  Q\right)  \text{Lip}\left(
\nu\right)  \leq11L\text{Lip}\left(  Q\right)  \text{Lip}\left(  \nu\right)
.
\]

\medskip

\noindent Using this, and property $\left(  2\right)  $ of $\delta_{n},$ we have,

\medskip%

\begin{equation}
\text{Lip}\left(  T_{n}\nu_{n}-f-\delta_{n}\right)  \leq11L\text{Lip}\left(
Q\right)  \text{Lip}\left(  \nu\right)  +L+C_{0}\varepsilon^{\prime}%
\leq13C_{0}L\text{Lip}\left(  Q\right)  \text{Lip}\left(  \nu\right)  .
\end{equation}

\medskip

\noindent Next, for $x\in D_{n},$ $\overline{\nu}_{n}\left(  x\right)  =1,$
and again by property $\left(  2\right)  $ of $\delta_{n},$
\begin{equation}
\text{Lip}\left(  \left(  T_{n}\overline{\nu}_{n}-f-\delta_{n}\right)
\mid_{D_{n}}\right)  =\ \text{Lip}\left(  \left(  T_{n}-f-\delta_{n}\right)
\mid_{D_{n}}\right)  \leq\varepsilon^{\prime}+C_{0}\varepsilon^{\prime}%
\leq2C_{0}\varepsilon^{\prime}.
\end{equation}
Next we compute, using $\left(  2.4\right)  ,$
\begin{align*}
\left\Vert \left(  T_{n}\left(  x\right)  \left(  \nu_{n}-\overline{\nu}%
_{n}\right)  \right)  ^{\prime}\right\Vert  &  =\left\Vert T_{n}^{\prime
}\left(  x\right)  \right\Vert \left\vert \nu_{n}\left(  x\right)
-\overline{\nu}_{n}\left(  x\right)  \right\vert +\left\vert T_{n}\left(
x\right)  \right\vert \left\Vert \nu_{n}^{\prime}\left(  x\right)
-\overline{\nu}_{n}^{\prime}\left(  x\right)  \right\Vert \\
& \\
&  \leq L\frac{\varepsilon^{\prime}r/2LL_{\varphi}}{1+\left\Vert
x-x_{n}\right\Vert }+\left(  L\left\Vert x-x_{n}\right\Vert +2\right)
\frac{\text{Lip}\left(  Q\right)  \varepsilon^{\prime}r/2LL_{\varphi}%
}{1+\left\Vert x-x_{n}\right\Vert }\\
& \\
&  \leq\varepsilon^{\prime}r/2+\text{Lip}\left(  Q\right)  \ \varepsilon
^{\prime}r/2+\text{Lip}\left(  Q\right)  \ \varepsilon^{\prime}r/L\\
& \\
&  \leq2\text{Lip}\left(  Q\right)  \varepsilon^{\prime}.
\end{align*}

\medskip

\noindent It follows from this and $\left(  2.18\right)  $ that,%

\begin{equation}
\text{Lip}\left(  \left(  T_{n}\nu_{n}-f-\delta_{n}\right)  \mid_{D_{n}%
}\right)  \leq2C_{0}\varepsilon^{\prime}+2\text{Lip}\left(  Q\right)
\varepsilon^{\prime}r\leq4C_{0}\text{Lip}\left(  Q\right)  \varepsilon
^{\prime}.
\end{equation}

\medskip

\noindent Finally we turn to $\left\Vert g^{\prime}\left(  x\right)
-f^{\prime}\left(  x\right)  \right\Vert $ with the help of the above
estimates. Again fix $x\in X,$ and we obtain,%

\begin{align*}
&  \left\Vert g^{\prime}\left(  x\right)  -f^{\prime}\left(  x\right)
\right\Vert \\
& \\
&  =\left(  \frac{\left\Vert \left\{  \left(  T_{n}\left(  x\right)  \nu
_{n}\left(  x\right)  -f\left(  x\right)  -\delta_{n}(x)\right)  \varphi
_{n}\left(  x\right)  \right\}  _{n=1}^{\infty}\right\Vert _{c_{0}}%
}{\left\Vert \left\{  \varphi_{n}\left(  x\right)  \right\}  _{n=1}^{\infty
}\right\Vert _{c_{0}}}\right)  ^{\prime}\\
& \\
&  =\frac{1}{\left\Vert \left\{  \varphi_{n}\left(  x\right)  \right\}
_{n=1}^{\infty}\right\Vert _{c_{0}}^{2}}\ \times\\
& \\
&  \left(  \left\Vert \left\{  \varphi_{n}\left(  x\right)  \right\}
_{n=1}^{\infty}\right\Vert _{c_{0}}\left\Vert \left\{  \left(  T_{n}\left(
x\right)  \nu_{n}\left(  x\right)  -f\left(  x\right)  -\delta_{n}(x)\right)
\varphi_{n}\left(  x\right)  \right\}  _{n=1}^{\infty}\right\Vert _{c_{0}%
}^{\prime}\right. \\
& \\
&  \left.  -\left\Vert \left\{  \varphi_{n}\left(  x\right)  \right\}
_{n=1}^{\infty}\right\Vert _{c_{0}}^{\prime}\left\Vert \left\{  \left(
T_{n}\left(  x\right)  \nu_{n}\left(  x\right)  -f\left(  x\right)
-\delta_{n}(x)\right)  \varphi_{n}\left(  x\right)  \right\}  _{n=1}^{\infty
}\right\Vert _{c_{0}}\right)
\end{align*}

\medskip

\noindent Let us first consider
\begin{align*}
&  \left(  \left(  T_{n}\left(  x\right)  \nu_{n}\left(  x\right)  -f\left(
x\right)  -\delta_{n}(x)\right)  \varphi_{n}\left(  x\right)  \right)
^{\prime}\\
& \\
&  =\left(  T_{n}\left(  x\right)  \nu_{n}\left(  x\right)  -f\left(
x\right)  -\delta_{n}(x)\right)  ^{\prime}\varphi_{n}\left(  x\right)
+\left(  T_{n}\left(  x\right)  \nu_{n}\left(  x\right)  -f\left(  x\right)
-\delta_{n}(x)\right)  \varphi_{n}^{\prime}\left(  x\right)  .
\end{align*}
For the first term, and $n\in\mathcal{N},$ we have, using property $\left(
2\right)  $ of Lemma \ref{lemma phi} and $\left(  2.19\right)  ,$%

\begin{align*}
\left\Vert \left(  T_{n}\left(  x\right)  \nu_{n}\left(  x\right)  -f\left(
x\right)  -\delta_{n}(x)\right)  ^{\prime}\varphi_{n}\left(  x\right)
\right\Vert  &  \leq\left\Vert \left(  T_{n}\left(  x\right)  \nu_{n}\left(
x\right)  -f\left(  x\right)  -\delta_{n}(x)\right)  ^{\prime}\right\Vert
\varphi_{n}\left(  x\right) \\
& \\
&  \leq4C_{0}\text{Lip}\left(  Q\right)  \varepsilon^{\prime}\left(
1+\varepsilon_{1}\right) \\
& \\
&  \leq8C_{0}\text{Lip}\left(  Q\right)  \varepsilon^{\prime}.
\end{align*}

\medskip

\noindent For $n\notin\mathcal{N},$ using Lemma \ref{lemma phi} $\left(
4^{\prime}\right)  $ and $\left(  2.17\right)  ,$ we obtain,%

\begin{align*}
\left\Vert \left(  T_{n}\left(  x\right)  \nu_{n}\left(  x\right)  -f\left(
x\right)  -\delta_{n}(x)\right)  ^{\prime}\varphi_{n}\left(  x\right)
\right\Vert  &  \leq\left\Vert \left(  T_{n}\left(  x\right)  \nu_{n}\left(
x\right)  -f\left(  x\right)  -\delta_{n}(x)\right)  ^{\prime}\right\Vert
\varphi_{n}\left(  x\right) \\
& \\
&  \leq13C_{0}L\text{Lip}\left(  Q\right)  \text{Lip}\left(  \nu\right)
\left(  \varepsilon^{\prime}/6C_{0}L\text{Lip}\left(  \nu\right)  \right) \\
& \\
&  \leq3\text{Lip}\left(  Q\right)  \varepsilon^{\prime}.
\end{align*}

\medskip

\noindent In any event, $\left\Vert \left(  T_{n}\left(  x\right)  \nu
_{n}\left(  x\right)  -f\left(  x\right)  -\delta_{n}(x)\right)  ^{\prime
}\varphi_{n}\left(  x\right)  \right\Vert _{c_{0}}\leq8C_{0}A_{1}$Lip$\left(
Q\right)  \varepsilon^{\prime}.$

\medskip

\noindent Next we consider the second term, $\left(  T_{n}\left(  x\right)
\nu_{n}\left(  x\right)  -f\left(  x\right)  -\delta_{n}(x)\right)
\varphi_{n}^{\prime}\left(  x\right)  .$ For $n\in\mathcal{N},$ from the
estimate giving $\left(  2.13\right)  ,$ we have,%
\begin{align*}
\left\Vert \left(  T_{n}\left(  x\right)  \nu_{n}\left(  x\right)  -f\left(
x\right)  -\delta_{n}(x)\right)  \varphi_{n}^{\prime}\left(  x\right)
\right\Vert  &  \leq\left\vert T_{n}\left(  x\right)  \nu_{n}\left(  x\right)
-f\left(  x\right)  -\delta_{n}(x)\right\vert \left\Vert \varphi_{n}^{\prime
}\left(  x\right)  \right\Vert \\
& \\
&  \leq6\varepsilon^{\prime}r/L_{\varphi}\left(  L_{\varphi}\right)
\leq6\varepsilon^{\prime}.
\end{align*}

\medskip

\noindent For $n\notin\mathcal{N},$ just as in $\left(  2.16\right)  $ we obtain,%

\begin{align*}
\left\Vert \left(  T_{n}\left(  x\right)  \nu_{n}\left(  x\right)  -f\left(
x\right)  -\delta_{n}(x)\right)  \varphi_{n}^{\prime}\left(  x\right)
\right\Vert  &  \leq\left\vert T_{n}\left(  x\right)  \nu_{n}\left(  x\right)
-f\left(  x\right)  -\delta_{n}(x)\right\vert \left\Vert \varphi_{n}^{\prime
}\left(  x\right)  \right\Vert \\
& \\
&  \leq\varepsilon^{\prime}.
\end{align*}

\noindent Hence, altogether we see that,
\begin{align*}
\left\Vert \left(  \left(  T_{n}\left(  x\right)  \nu_{n}\left(  x\right)
-f\left(  x\right)  -\delta_{n}(x)\right)  \varphi_{n}\left(  x\right)
\right)  ^{\prime}\right\Vert _{c_{0}}  &  \leq8C_{0}A_{1}\text{Lip}\left(
Q\right)  \varepsilon^{\prime}+6A_{1}\varepsilon^{\prime}\\
& \\
&  \leq14C_{0}A_{1}\text{Lip}\left(  Q\right)  \varepsilon^{\prime}.
\end{align*}

\medskip

\noindent Lastly, we examine $\left\Vert \left\{  \varphi_{n}\left(  x\right)
\right\}  _{n=1}^{\infty}\right\Vert _{c_{0}}^{\prime}\left\Vert \left\{
\left(  T_{n}\left(  x\right)  \nu_{n}\left(  x\right)  -f\left(  x\right)
-\delta_{n}(x)\right)  \varphi_{n}\left(  x\right)  \right\}  _{n=1}^{\infty
}\right\Vert _{c_{0}}.$ Recall our estimate of $\left\vert f-g\right\vert $
found $\left\Vert \left\{  \left(  T_{n}\left(  x\right)  \nu_{n}\left(
x\right)  -f\left(  x\right)  -\delta_{n}(x)\right)  \varphi_{n}\left(
x\right)  \right\}  _{n=1}^{\infty}\right\Vert _{c_{0}}\leq3A_{1}%
\varepsilon^{\prime}r.$ Therefore we have,%

\begin{align*}
&  \left\Vert \left\{  \varphi_{n}\left(  x\right)  \right\}  _{n=1}^{\infty
}\right\Vert _{c_{0}}^{\prime}\left\Vert \left\{  \left(  T_{n}\left(
x\right)  \nu_{n}\left(  x\right)  -f\left(  x\right)  -\delta_{n}(x)\right)
\varphi_{n}\left(  x\right)  \right\}  _{n=1}^{\infty}\right\Vert _{c_{0}}\\
& \\
&  \leq A_{1}L_{\varphi}\left\Vert \left\{  \left(  T_{n}\left(  x\right)
\nu_{n}\left(  x\right)  -f\left(  x\right)  -\delta_{n}(x)\right)
\varphi_{n}\left(  x\right)  \right\}  _{n=1}^{\infty}\right\Vert _{\infty}\\
& \\
&  \leq A_{1}\frac{L_{1}\text{Lip}\left(  Q\right)  }{r}\left(  5A_{1}%
\varepsilon^{\prime}r\right)  =5A_{1}^{2}L_{1}\text{Lip}\left(  Q\right)
\varepsilon^{\prime}%
\end{align*}

\medskip

\noindent Finally, because $\left\Vert \left\{  \varphi_{n}\left(  x\right)
\right\}  _{n=1}^{\infty}\right\Vert _{c_{0}}^{2}\geq\left\Vert \left\{
\varphi_{n}\left(  x\right)  \right\}  _{n=1}^{\infty}\right\Vert _{\infty
}^{2}\geq1/4$ as noted above, putting all the above estimates together
yields,
\begin{align*}
\left\Vert g^{\prime}\left(  x\right)  -f^{\prime}\left(  x\right)
\right\Vert  &  \leq\frac{\left(  2A_{1}\right)  14C_{0}A_{1}\text{Lip}\left(
Q\right)  \varepsilon^{\prime}+5A_{1}^{2}L_{1}\text{Lip}\left(  Q\right)
\varepsilon^{\prime}}{1/4}\\
& \\
&  \leq\left(  132C_{0}A_{1}^{2}L_{1}\text{Lip}\left(  Q\right)  \right)
\varepsilon^{\prime}<\varepsilon.\ \ \square
\end{align*}

\end{document}